\documentclass[12pt]{article}
\usepackage[cp866]{inputenc}
\usepackage{citehack}
\usepackage{amsmath}
\usepackage{amsfonts}
\usepackage{amssymb}

\def\Simil{\text{\rm Sim}}
\def\End{\text{\rm End}}
\def\spa{\text{\rm span}}
\def\Hom{\text{\rm Hom}}
\def\rk{\text{\rm rk}}

\def\Real{\mathbb{R}}
\def\g{\mathfrak{g}}
\def\h{\mathfrak{h}}
\def\so{\mathfrak{so}}

\def\gl{\mathfrak{gl}}
\def\su{\mathfrak{su}}
\def\u{\mathfrak{u}}
\def\sp{\mathfrak{sp}}
\def\f{\mathfrak{f}}
\def\t{\mathfrak{t}}
\def\z{\mathfrak{z}}
\def\R{\mathcal{R}}
\def\P{\mathcal{P}}

\begin{document}
\title{The spaces of curvature tensors for  holonomy
algebras of Lorentzian manifolds}

\author{Anton S. Galaev \thanks{EMail: \tt  galaev@mathematik.hu-berlin.de}\\
Institut f\"ur Mathematik, Humboldt-Universit\"at zu Berlin\\ Sitz: Rudower Chaussee 25, 10099, Berlin, Germany}

\maketitle\vskip-50ex
 {\renewcommand{\abstractname}{Abstract}\begin{abstract}
  The holonomy algebra $\g$ of an indecomposable Lorentzian
(n+2)-dimensional
manifold $M$  is  a weakly-irreducible subalgebra of the Lorentzian algebra
$\so_{1,n+1}$.
  L. Berard Bergery and A. Ikemakhen  divided weakly-irreducible
  not irreducible
subalgebras  into 4 types and associated
with each such subalgebra $\g$ a subalgebra $\h\subset \so_n$
of the orthogonal Lie algebra.
We give a description of the spaces $\R(\g)$ of the curvature
tensors for algebras of each type
in terms of the space $\P(\h)$ of $\h$-valued 1-forms
on $\Real^n$ that satisfy the Bianchi identity and reduce the
classification
of the holonomy algebras of Lorentzian manifolds to the classification
of irreducible subalgebras $\h$ of $so(n)$
with $L(\P(\h))=\h$.
  We prove that for
$n\leq 9$ any such subalgebra $\h$ is  the holonomy algebra of a
Riemannian manifold.
This gives a classification of the holonomy algebras for Lorentzian
manifolds $M$ of dimension $\leq 11$.
\end{abstract}


{\bf Mathematical subject codes:} 53c29, 53c50

\section*{Introduction}

The connected irreducible holonomy groups of pseudo-Riemannian
manifolds have been classified by M. Berger,
see  
\cite{Ber}. The classification problem for not irreducible holonomy groups
 is still open.
The main difficulty is that the holonomy group can preserve
an isotropic subspace of the tangent space.

We tackle the  classification problem for
the holonomy algebras of Lorentzian manifolds,
i.e. the Lie algebras of the holonomy groups.
There are some partial results in this direction (see 
\cite{B-I,Ik,Le1,Le2}).
In \cite{As} the classification of holonomy algebras for
4-dimensional Lorentzian manifolds was given.

Wu's theorem (see 
\cite{Wu}) reduces the classification problem
for holonomy algebras to the classification of weakly-irreducible holonomy
algebras
(i.e. algebras that preserve no nondegenerate proper
subspace of the tangent space).

If a holonomy algebra is irreducible, then it is weakly-irreducible.
The Berger list (\cite{Ber})
of irreducible holonomy algebras of
pseudo-Riemannian manifolds  shows that the only irreducible
holonomy algebra of Lorentzian manifolds is
$\so(1,m-1)$, see also 
\cite{D-O}.

We study weakly-irreducible holonomy algebras that are not
irreducible.

Let $(V,\eta)$ be an $(n+2)$-dimensional Minkowski vector space,
where
$\eta$ is a metric of signature $(1,n+1)$.
Using $\eta$ we  identify the space $V$ with the dual space $V^*$.
Then the Lorentzian algebra $\so(V)$
is identified with the space $V \wedge V$ of bivectors.
 Denote by $\so(V)_{\Real p}$
the subalgebra of $\so(V)$ that preserves an isotropic
line $\Real p$, where $p\in V$. Denote by $E$ a vector subspace
$E\subset V$ such that $(\Real p)^{\bot_\eta}=\Real p\oplus E$.
The vector space $E$ is an Euclidean space with respect to the
inner product $-\eta|_E.$ Denote by $q$ an isotropic vector $q\in
V$ such that $\eta(q,E)=0$ and $\eta(p,q)=1$.
We have $$\so(V)_{\Real p}=\Real p\wedge q+p\wedge E+\so(E).$$
Any weakly-irreducible and not irreducible subalgebra of $\so(V)$
 is conjugated to a subalgebra $\g$ of $\so(V)_{\Real p}$. We denote
 by $\h^{\g}$ the projection of such subalgebra
 $\g$ to $\so(E)$ with respect to the above decomposition
  and call $\h^{\g}$ the
 { \it orthogonal part} of the Lie algebra $\g$.

Conversely, for any subalgebra $\h\subset\so(E)$ we construct two
Lie algebras
$$\g^\h_1=\Real p\wedge q +\h+p\wedge E$$
and $$\g^\h_2=\h+p\wedge E$$
with the orthogonal part $\h$.
Moreover, if the center $\z(\h)$ of $\h$ is non-trivial, then
any non-zero linear map $$\varphi:\z(\h)\to\Real$$ defines the Lie algebra
$$\g_3^{\h ,\varphi}= p\wedge E + \{A + \varphi(A)p\wedge q:\, A \in \h \}$$
with the orthogonal part $\h$. Here $\varphi$ is considered as
the linear map $\varphi :\h=\z(\h)\oplus\h'\to\Real$
that vanishes on the commutant $\h'$ of $\h$.

  Suppose moreover that the subalgebra $\h \subset
   \so(E)$ acts trivially on a subspace $E_0 \neq \{0\}$, such that we
   can consider $\h$ as a subalgebra of $\so(E_1)$, where
   $E = E_0\oplus E_1$ is the orthogonal decomposition. Then any
   surjective linear map
   $$ \psi : \z (\h) \to E_0 $$
   extended to $\h$ by $\psi(\h')=0$ defines the Lie
   algebra
   $$\g_4^{\h,\psi}= p \wedge E_1 + \{A + p\wedge \psi(A):\,A \in \h\}$$
     with the orthogonal part $\h$.
We call Lie algebras
 $\g^\h_1$, $\g^\h_2$, $\g^{\h,\varphi}_3$ and $\g^{\h,\psi}_4$
 with the orthogonal part $\h \subset\so(E)$
 the algebras of type 1,2,3 and 4 respectively.

These Lie algebras were considered by L. Berard Bergery and
A. Ikemakhen 
\cite{B-I}, who proved that {\it the Lie algebras of the form
$\g_1^\h$, $\g_2^\h$, $\g_3^{\h,\varphi}$, and $\g_4^{\h,\psi}$ exhaust
all weakly-irreducible subalgebras of $\so(V)_{\Real p}$} (theorem 1).
The other result is that {\it the orthogonal part of the holonomy algebra of a
Lorentzian manifold satisfies a Borel-Lichnerowicz-type decomposition property} (theorem 2).

{\bf Remark.} Note that the Lie algebra $\so(V)_{\Real p}=\Real p\wedge q+p\wedge E+\so(E)$
is isomorphic to the tangent Lie algebra for the Lie group $\Simil E$ of similarity transformations of $E$,
the elements $\lambda p\wedge q$ and $p\wedge u$ correspond to the homothetic transformation $v\mapsto \lambda v$
and to the shift $v\mapsto v+u$ respectively, here $\lambda\in\Real$ and $u,v\in E$. In another paper we will give a geometrical
interpretation to the result of L. Berard Bergery and A. Ikemakhen.

Let $\g\subset\so(V)_{\Real p}$ be a subalgebra.
Recall that the space of curvature tensors of type $\g$ is defined as
         the space $\R(\g)$ of $\g$-valued 2-forms
         on $V$ that satisfy the Bianchi identity. We
         denote by
 $$L(\R(\g))=\text{span}(\{R(u\wedge v):\, R \in R(\g),\, u,v \in V\})$$
         the vector subspace of $\g$ spanned by curvature operators
         from $\R(\g)$. If $\g$ is the holonomy
         algebra of an indecomposable Lorentzian manifold, then
         $$L(\R (\g))=\g. \,\,\,\,\,\,\,\,\,\,\,\,\,\,\,\,\,\,\,\,\, (*)$$
          A weakly-irreducible subalgebra $\g \subset \so(V)_{\Real p}$
          that satisfies
          (*)  is called a {\it  Berger algebra}.

In this  paper we give a description of the spaces $\R(\g)$ of curvature
tensors for  weakly-irreducible
 subalgebras $\g \subset \so(V)_{\Real p}$ of each type in terms of
the orthogonal part $\h^{\g} \subset \so(E)$
and we reduce the classification
of Berger algebras  to the classification of irreducible subalgebras
$\h\subset\so(E)$ that satisfy some conditions (weak-Berger algebras).

 More
          precisely,  for any subalgebra $\h\subset \so(E)$
          we define the space
          $$\P(\h) =\{P\in \text{Hom}(E,\h):
\,\eta(P(u)v,w)+\eta(P(v)w,u)+\eta(P(w)u,v)=0 \text{ for all }u,v,w\in E\}$$
          of $\h$-valued 1-forms on $E$ that satisfy the
          Bianchi identity  and denote by
          $$L(\P(\h) ) = \text{span}(\{P(u):P\in\P(\h),\,u\in E\})$$
          the vector subspace of $\h$  spanned by tensors  $P\in\P(\h)$.
          We call $\P(\h)$ the {\it space of weak-curvature tensors of type}
          $\h$.

A subalgebra $\h \subset \so(E)$ is called a {\it weak-Berger algebra}
          if $L(\P(\h))=\h$.

We give a description of the
spaces of curvature tensors $\R(\g)$ for
algebras of each type associated with a
given orthogonal part $\h \subset \so(E)$
in terms of the space $\P(\h)$ of
weak-curvature tensors (theorem 3).

Corollary 1 shows that {\it a weakly-irreducible
subalgebra $\g\subset\so(V)_{\Real p}$
is a Berger algebra iff
$\h^\g$ is a weak-Berger algebra.}

Note that the direct sum $\h_1\oplus\h_2\subset \so(E_1)\oplus\so(E_2)$
of two weak-Berger algebras  $\h_i\subset\so(E_i), i=1,2 $
is a weak-Berger algebra.

In part (I) of theorem 4 we prove that\\ {\it if a subalgebra
$\h\subset\so(E)$ is a weak-Berger algebra, then
there exists an orthogonal decomposition $E=E_0\oplus
E_1\oplus\cdots\oplus E_r$ and the corresponding decomposition into
the direct sum of ideals
$\h=\{0\}\oplus\h_1\oplus\h_2
\oplus\cdots\oplus\h_r$ such that
$\h_i(E_j)=0$ if $i\neq j,$
$\h_i\subset\so(E_i),$ and $\h_i$ acts
irreducibly on $E_i.$
}
This result is stronger than theorem 2,
without this it was necessary
to suppose that a Berger algebra satisfies the conclusion of
theorem 2.
My attention to this statement  was taken
by A.J. Di Scala.

Part (II) of theorem 4 states that {\it if $\h\subset\so(E)$ is
a weak-Berger algebra, then we have
$$\P(\h)=\P(\h_1)\oplus\cdots\oplus\P(\h_r).$$} This
implies that {\it  $\h_i$ is a
weak-Berger algebra for $i=1,...,k$.}

Thus the classification problem for the Berger algebras
$\g\subset\so(V)_{\Real p}$
 is reduced to
 the classification of irreducible weak-Berger algebras $\h\subset\so(E)$.

We prove that {\it  if $\h\subset\so(E)$ is the holonomy algebra
of a Riemannian  manifold, then it is a weak-Berger
algebra.}

Using theory of representation of compact Lie algebras,
we prove the converse statement in the case
when $\dim E \leq 9$. Theorem 5 states that
{\it if $\dim V\leq 11,$ then
a weakly-irreducible subalgebra $\g\subset\so(V)_{\Real p}$
is a Berger algebra iff
the algebra $\h^\g$ is
the holonomy
algebra of a Riemannian manifold}.
This gives a classification of
Berger algebras for Lorentzian manifolds of dimension $\leq 11$,
which can be stated in the following way.

Let $n_0,n_1,...,n_r$ be positive integers such that
$2\leq n_1\leq\cdots\leq n_r$ and $n_0+n_1+\cdots+n_r=n$.
Let $\h_i\subset \so_{n_i}$ be the holonomy algebra of an irreducible Riemannian manifold
 ($i=1,...,r$). The Lie algebras of the form
 $\h=\h_1\oplus\cdots\oplus\h_r$ exhaust
all weak-Berger subalgebras of $\so_n$. The Lie algebras of the
form
$\g_1^\h$, $\g_2^\h$, $\g_3^{\h,\varphi}$ and
$\g_4^{\h,\psi}$ (if $\g_3^{\h,\varphi}$ and $\g_4^{\h,\psi}$ exist)
exhaust all Berger algebras for Lorentzian manifolds of dimension
$n+2$. Note that for each $n>1$ there exists infinite number of
weakly-irreducible Berger subalgebras of $\so(V)_{\Real p}$. 

The full list of irreducible
holonomy algebras of Riemannian
manifolds of dimension $\leq 9$
is given in the table below.
In the table $\otimes$
stands for the tensor product of representations;
$\underline{\otimes}$
stands for the highest irreducible component of the corresponding
product.

\vskip 0.3cm

\begin{tabular}{|c|c|}
\hline
 $n$ & Irreducible weak-Berger subalgebras of $\so_n$\\ \hline
  $n=1$ &    \\ \hline
  $n=2$ & $\so_2$   \\ \hline
  $n=3$ & $\so_3$  \\ \hline
  $n=4$ & $\so_4$,
  $\su_2$,
  $\u_2$

    \\ \hline
  $n=5$ & $\so_5$,
  $\underline{\otimes}^2\so_3$  \\ \hline
  $n=6$ & $\so_6$,
  $\su_3,$
   $\u_3$
  \\ \hline
  $n=7$ & $\so_7$, $\g_2$
     \\ \hline

   $n=8$ & $\so_8$,
    $\su_4$,
    $\u_4$,

$\sp_2$,
    $\sp_2\otimes\sp_1$,
    $\underline{\otimes}^3\so_3\otimes\so_3$,

    $\underline{\otimes}^2\su_3$

    \\      \hline

$n=9$ & $\so_9$,
$\so_3\otimes\so_3$
 \\      \hline
  \end{tabular}

\vskip 0.4cm

Recall that the holonomy group of an indecomposable Lorentzian
manifold can be not closed. In 
\cite{B-I} it was shown that
the connected Lie subgroups of $SO_{1,n+1}$ corresponding to Lie
algebras of type 1 and 2 are closed; the connected Lie subgroup
of $SO_{1,n+1}$ corresponding to a Lie
algebra of type 3 (resp. 4) is closed if and only if the connected
subgroup of $\so_n$ corresponding to the subalgebra
$\ker\varphi\subset\z(\h)$ (resp. $\ker\psi\subset\z(\h)$) is
closed. We give a criteria for Lie groups corresponding to Lie
algebras of type 3 and 4 to be closed in terms of the Lie algebras
$\ker\varphi$ and $\ker\psi$.

Let $\h\subset\so_n$ be a weak-Berger algebra such that $\z(\h)\neq\{0\}$.
We have $\h=\h_1\oplus\cdots\oplus\h_r$, where $\h_i$ are
irreducible weak-Berger algebras.
We see that $\z(\h)=\z(\h_1)\oplus\cdots\oplus\z(\h_r)$ and
$\dim\z(\h_i)=0$ or $1$, $i=1,...,r$. Hence we can identify
$\z(\h)$ with $R^m$, where $m=\dim\z(\h)$.

We prove that\\ {\it the connected Lie group corresponding
to a Lie algebra of type 3 (resp. 4) is closed if
and only if there exists a basis $v_1,...,v_l$
($l=\dim\ker\varphi$ or $\dim\ker\psi$)
of the vector space $\ker\varphi$ (resp. $\ker\psi$)
such that the coordinates of the vector $v_i$ with respect to the canonical
basis of $\Real^m$ are integer for $i=1,...,l$.}

{\bf Remark.}  Recently Thomas Leistner proved that {\it a subalgebra
$\h\subset\so(E)$ is a weak-Berger algebra iff $\h$ is the
holonomy algebra of a Riemannian manifold}, see 
\cite{Le1,Le2,Le3}.

From this result and corollary 1 it follows that
{\it a weakly-irreducible
subalgebra $\g\subset\so(V)_{\Real p}$
is a Berger algebra iff
$\h^\g$ is the holonomy algebra of a Riemannian manifold.}

\centerline{\bf Acknowledgments}
I am grateful to D.V. Alekseevsky for
introducing me to holonomy algebras of Lorentzian manifolds and
for his helpful suggestions. I also thank my
teacher M.V. Losik for long hours of discussions
and help in preparation of this paper.

\section{Preliminaries}

Let $(V,\eta)$ be a Minkowski space of dimension $n+2$,
where  $\eta$ is a metric on $V$ of
signature $(1,n+1)$. We fix a basis
$p,e_1,...,e_n,q$ of $V$ such that the Gram matrix
of $\eta$ has the form
$\left (\begin{array}{ccc}
0 & 0 & 1\\ 0 & -E_n & 0 \\ 1 & 0 & 0 \\
\end{array}\right)$,
 where $E_n$ is the $n$-dimensional identity matrix.

Let $E\subset V$ be the vector subspace spanned by $e_1,...,e_n$.
We will consider $E$ as an Euclidean space with the metric
$-\eta|_E.$

Denote by $\so(V)$ the Lie algebra of all $\eta$-skew symmetric
endomorphisms of $V$ and by $\so(V)_{\Real p}$ the subalgebra
of $\so(V)$ that preserves the line $\Real p$.

The algebra $\so(V)_{\Real p}$ can be identified with
the following matrix algebra:
  $$\so(V)_{\Real p}=\left\{ \left (\begin{array}{ccc}
a & X^t & 0\\ 0 & A & X \\ 0 & 0 & -a \\
\end{array}\right):\, a\in \Real,\, X\in \Real^n,\,A \in \so(n)
 \right\}.$$

We identify the dual vector space
$V^*$ with $V$ using $\eta$.
Hence we can identify $\End V=V\otimes V^*$ with $V\otimes V$.
In particular, we identify
$\so(V)$ with $V\wedge V=
\text{span}(\{u\wedge v=u\otimes v-v\otimes u:u,v\in V\}).$
Similarly, we identify $\so(E)$ with $E\wedge E$ and consider
$\so(E)$ as a subspace of $\so(V)$ that acts trivially on $\Real p\oplus\Real q$.

For $a\in \Real$, the endomorphism
$ap\wedge q$ has the matrix
$\left (\begin{array}{ccc}
a & 0 & 0\\ 0 & 0 & 0 \\ 0 & 0 & -a \\
\end{array}\right)\in \so(V)_{\Real p};$ for $X\in E$, the
endomorphism $p\wedge X$ has the matrix
$\left (\begin{array}{ccc}
0 & X^t & 0\\ 0 & 0 & X \\ 0 & 0 & 0 \\
\end{array}\right)\in \so(V)_{\Real p}.$ Thus we see that
$$\so(V)_{\Real p}=E\wedge E+p\wedge E+\Real p\wedge q \text{ is
a direct sum of the subalgebras.}$$

{\bf Definition 1.}
{\it A subalgebra $\g\subset \so(V)$ is
called irreducible if it preserves no proper
subspace of $V$; $\g$ is called weakly-irreducible if
it preserves no nondegenerate proper subspace of $V$.
}

Obviously, if $\g\subset \so(V)$ is
irreducible, then it is weakly-irreducible.
If $\g\subset \so(V)$ preserves
a degenerate proper subspace $U\subset V$, then it preserves the
isotropic line $U\bigcap U^\bot$;
any such algebra is conjugated to a subalgebra of $\so(V)_{\Real p}$.

{\bf Definition 2.}
{\it Let $W$ be a vector space and
$\f\subset \gl(W)$ a subalgebra.
Put}
$$\R(\f)=\{R\in\text{Hom}(W\wedge W,\f):
R(u\wedge v)w+R(v\wedge w)u+R(w\wedge u)v=0 \text{ }\text{\it for all } u,v,w\in W\}.$$
{\it The set $\R(\f)$ is called the space of
curvature tensors of type $\f.$
Denote by $L(\R(\f))$ the vector
subspace of $\f$ spanned by $R(u\wedge v)$ for all
$R\in\R(\f),$ $u,v\in W,$}
$$L(\R(\f))=\spa(\{R(u\wedge v):R\in\R(\f),\,u,v\in W\}).$$

Let $\g\subset \so(V)$ be a subalgebra.
Recall that a curvature tensor $R\in\R(\g)$ satisfies
the following property
\begin{equation}
 \eta (R(u\wedge v)z,w)=\eta
(R(z\wedge w)u,v)\, \text{ for all } u,v,z,w\in V.\label{property}
\end{equation}

Let $(M,g)$ be a Lorentzian manifold of
                    dimension $n+2$  and $\g$ the holonomy
                    algebra (that is the Lie algebra of the holonomy group)
                    at a point $x$. By Wu's theorem (see \cite{Wu}) $(M,g)$ is
                    locally
                    indecomposable, i.e. is not locally a product
                    of two pseudo-Riemannian manifolds if and only
                    if the holonomy algebra  $\g$ is weakly-irreducible. If the holonomy algebra $\g$
                    is irreducible , then $\g =\so(1,n+1)$.
                    So we may assume that it is reducible and weakly-irreducible. Then it
                    preserves an isotropic line $\ell \subset
                    T_xM$. We can identify the tangent space
                    $T_xM$ with $V$ such that $\ell$ corresponds
                    to the line $\Real p$. Then $\g$ is
                    identified with  weakly-irreducible
                    subalgebra of  $\so(V)_{\Real p}$.

                     We need the following

{\bf Proposition 1.}
{\it Let $\g$ be the holonomy algebra of a Lorentzian manifold.
 Then $$L(\R(\g))=\g.$$}

{\bf Proof.}
The inclusion $L(\R(\g))\subset\g$ is obvious.

Let $R$ be the curvature tensor of $(M,g).$
Theorem of Ambrose and Singer
(\cite{Am-Si}) states that the vector space $\g$ is generated by all  endomorphisms
$$(\tau(\lambda))^{-1}\circ R_{\lambda(b)}(\tau (\lambda)(X), \tau (\lambda)(Y))
\circ \tau (\lambda):\,T_xM\to T_xM,$$ where
$\lambda:\,[a,b]\to M$ is a piecewise smooth curve in $M$ such that $\lambda(a)=x$,
$\tau(\lambda)$ is the parallel
transport along $\lambda$, and $X,Y\in T_{\lambda(a)}M$.
Obviously, the above transformations  are curvature tensors of type
$\mathfrak{g},$
hence, $\mathfrak{g}\subset L(\mathcal{R}(\mathfrak{g})).$
Thus, $\mathfrak{g}=L(\mathcal{R}(\mathfrak{g})).$ $\Box$

{\bf Definition 3.}
{\it A weakly-irreducible subalgebra $\g\subset \so(V)_{\Real p}$
is called a
Berger algebra if
$L(\R(\g))=\g$.
}

{\bf Definition 4.}
{\it Let $\h\subset\so(E)$ be a subalgebra.  Put
$$\P(\h)=\{P\in \Hom (E,\h):
\,\eta(P(u)v,w)+\eta(P(v)w,u)+\eta(P(w)u,v)=0\text{ }\text{\it for all }u,v,w\in
E\}.$$
We call $\P(\h)$ the space of weak-curvature tensors of type
$\h$.
A subalgebra $\h\subset\so(E)$ is called a weak-Berger algebra if
$L(\P(\h))=\h$, where
$$L(\P(\h))=\spa(\{P(u):P\in\P(\h),\,u\in E\})$$
is the vector subspace of $\h$ spanned by
$P(u)$ for all $P\in\P(\h)$ and $u\in E$.}

Let $\h\subset \so(E)$ be a subalgebra. Since
$\h$ is a compact Lie algebra, we have
$\h=\h'\oplus\z(\h)$ (the
direct sum of ideals), where
$\h'$ is the commutant of $\h$ and
$\z(\h)$
is the center of $\h$.

Consider a weakly-irreducible subalgebra $\g\subset \so(V)_{\Real p}$.
Let $\h^\g$ be the projection of $\g$ to $\so(E)$ with respect
to the decomposition $\so(V)_{\Real p}=\so(E)+p\wedge E+\Real p\wedge q$.

{\bf Definition 5.}
{\it The Lie algebra $\h^\g$ is called the orthogonal
part of $\g$.}

Conversely, with any subalgebra $\h\subset\so(E)$ we associate two
Lie algebras
$$\g^\h_1=\Real p\wedge q +\h+p\wedge E$$
and $$\g^\h_2=\h+p\wedge E.$$

Moreover, suppose $\z(\h)\neq\{0\}$.
Let $$\varphi:\,\z(\h) \to \Real$$ be a non-zero linear
map. Extend $\varphi$ to the linear map
$\varphi:\,\h \to \Real$ by putting $\varphi|_{\h'}=0.$
Then $$\g_3^{\h,\varphi}=p\wedge E+\{A+\varphi(A)p\wedge q:\,A\in
\h\}$$ is a Lie algebra with the orthogonal part $\h$.

Suppose moreover that we have an orthogonal decomposition
$E=E_0\oplus E_1$ such that\\ $E_0\neq\{0\}$, $\h\subset\so(E_1)$, and
$\dim\z(\h)\geq\dim E_0$.
Let $$\psi:\,\z(\h) \to E_0$$ be a surjective linear map.
As above, we extend $\psi$ to a linear map
$\psi:\h\to E_0$ by
 putting $\psi|_{\h'}=0.$
Then $$\g_4^{\h,\psi}=p\wedge E_1+\{A+p\wedge \psi(A):\,A\in \h\}$$
is a Lie algebra with the orthogonal part $\h$.

We call the Lie algebras
 $\g^\h_1$, $\g^\h_2$, $\g^{\h,\varphi}_3$ and $\g^{\h,\psi}_4$
the Lie algebras of type 1,2,3 and 4 respectively.

These Lie algebras were considered by L. Berard Bergery and
A. Ikemakhen, who proved the following fundamental results (see 
\cite{B-I}).

{\bf Theorem 1.}
{Let $\h\subset\so(E)$ be a subalgebra. Then (if $\g_3^{\h,\varphi}$ and $\g_4^{\h,\psi}$ exist) the subalgebras
$\g_1^\h$, $\g_2^\h$, $\g_3^{\h,\varphi}$,
$\g_4^{\h,\psi}\subset\so(V)_{\Real p}$ are weakly-irreducible.
Moreover, Lie algebras of the form $\g_1^\h$, $\g_2^\h$, $\g_3^{\h,\varphi}$ and
$\g_4^{\h,\psi}$ exhaust all weakly-irreducible subalgebras of
$\so(V)_{\Real p}$.}

{\bf Theorem 2.} {\it  Let $\g$ be the holonomy algebra of
a Lorentzian manifold.
Then there exists an orthogonal decomposition $E=E_0\oplus
E_1\oplus\cdots\oplus E_r$ and the corresponding decomposition into
the direct sum of ideals
$\h^\g=\{0\}\oplus\h_1\oplus\h_2
\oplus\cdots\oplus\h_r$ such that
$\h_i(E_j)=0$ if $i\neq j,$
$\h_i\subset\so(E_i),$ and $\h_i$ acts
irreducibly on $E_i.$
}

The metric $\eta$ on $V$ induces the metrics on $V\otimes V$ and
$V\wedge V$. Denote those metrics by $\eta\otimes\eta$ and $\eta\wedge\eta$
respectively.

Let $\theta:\,V\to V $ be an endomorphism, $u,v\in V.$
Then $\eta(\theta(u),v)=\eta\otimes\eta(\theta,u\otimes v).$
In particular, for $\theta\in V\wedge V$ we have
\begin{equation}
\eta(\theta(u),v)=1/2\eta\wedge\eta(\theta,u\wedge v).\label{property2}
\end{equation}

Let $R\in\R(\g).$ Combining \eqref{property}
and \eqref{property2}, we see that
$\eta\wedge\eta(R(u\wedge v),z\wedge w)=\eta\wedge\eta(R(z\wedge w),u\wedge v)$
for all $u,v,z,w\in V$.
This means
that the linear map $R:\,V\wedge V\to \g\subset V\wedge V$
is $\eta\wedge\eta$-symmetric.

Let $(E_1,\mu_1)$ and $(E_2,\mu_2)$ be two
Euclidean spaces. Let $f:\,E_1 \to E_2$ be a linear map.
Denote by $f^*:E_2\to E_1$ the dual linear map for $f$.
We identify the symmetric square $S^2(E)$ of $E$
with the space of all $\eta$-symmetric endomorphisms of $E$.

\section{Main results}

Let $\h\subset\so(E)$ be a subalgebra.
We will define some sets of endomorphisms,
in theorem 3 we will see that those sets consist of the
curvature tensors for appropriate algebras.

For any $\lambda\in \Real,$ $L\in\Hom(E,\Real),$ $T\in S^2(E)$ and
$P\in\P(\h)$ we define the endomorphisms
$$R^\lambda\in\Hom(V\wedge V,\g^\h_1),\, R^L\in\Hom(V\wedge V,\g^\h_1),\,
R^T\in\Hom(V\wedge V,p\wedge E) \text { and } R^P\in\Hom(V\wedge V,\g^\h_2)$$
by conditions
$$R^\lambda(p\wedge q)=\lambda p\wedge q,\,
R^\lambda|_{p\wedge E+q\wedge E+E\wedge E}=0,$$
$$R^L(q\wedge \cdot)=L(\cdot)p\wedge q,\,
R^L(p\wedge q)=p\wedge L^*(1),\,R^L|_{p\wedge E+E\wedge E}=0,$$
$$R^T(q\wedge \cdot)=p\wedge T(\cdot),\,
R^T|_{\Real p\wedge q+p\wedge E+E\wedge E}=0$$
and
$$R^{P}(q\wedge \cdot)=P(\cdot),\,
R^{P}|_{E\wedge E}=-1/2p\wedge P^*,\
R^{P}|_{\Real p\wedge q+p\wedge E}=0,$$
and define by
$$\R(\Real,\Real),\, \R(E,\Real),\,\R(E,E) \text{ and } \R(E,\h)$$
respectively the vector spaces of all such endomorphisms.

We have the isomorphisms
$\R(\Real,\Real)\simeq \Real,$ $\R(E,\Real)\simeq E,$ $\R(E,E)\simeq S^2(E)$
and $\R(E,\h)\simeq \P(\h)$.

Moreover, if a Lie algebra $\g^{\h,\varphi}_3$ exists, then for any
$P\in\P(\h)$  we define the endomorphism
$R^P\in\text{Hom}(V\wedge V,\g^{\h,\varphi}_3)$ by conditions
$$R^{P}(q\wedge \cdot)=P(\cdot)+\varphi(P(\cdot))p\wedge q,\,
R^{P}|_{E\wedge E}=-1/2p\wedge P^*,$$
$$R^{P}(p\wedge q)=-1/2p\wedge P^*(\varphi^*(1)),\,
R^{P}|_{p\wedge E}=0$$
and denote by $\R(E,\h,\varphi)$
the vector space of all such endomorphisms.

If a Lie algebra $\g^{\h,\psi}_4$ exists,
then for any
$P\in\P(\h)$  we define the endomorphism
$R^P\in\text{Hom}(V\wedge V,\g^{\h,\psi}_4)$ by conditions
$$R^{P}(q\wedge u_1)=P(u_1)+p\wedge\psi(P(u_1)) \text{ for all } u_1\in
E_1,\,R^P|_{E\wedge E}=-1/2p\wedge P^*,$$
$$R^{P}(q\wedge u_0)=-1/2p\wedge P^*(\psi^*(u_0))\text{ for all } u_0\in
E_0,\,R^{P}|_{\Real p\wedge q+p\wedge E}=0$$
and denote by $\R(E_1,\h,\psi)$
vector space of all such endomorphisms.

We have the isomorphisms $\R(E,\h,\varphi)\simeq\P(\h)$ and
$\R(E_1,\h,\psi)\simeq\P(\h)$.

Let $\g\subset\so(V)_{\Real p}$ be a weakly-irreducible
subalgebra, $\h^\g$ be the orthogonal part of $\g$ and
$\h\subset\h^\g$ be a subalgebra. Suppose $R\in\R(\h)$.
Extend the linear map $R:E\wedge E\to\h$ to the linear map
$R:V\wedge V\to\h$ by putting
$R^{P}|_{\Real p\wedge q+p\wedge E+q\wedge E}=0$. It is obvious
that $R\in\R(\g)$. We can write $\R(\h)\subset\R(\g)$.

{\bf Theorem 3.} {\it Let $\h\subset\so(E)$ be a subalgebra. Then
we have}\\
(I) $\R(\g_1^\h)=\R(\g_2^\h)\oplus
\R(E,\Real)\oplus\R(\Real,\Real);$\\
(II) $\R(\g_2^\h)=\R(\h)\oplus \R(E,\h)\oplus\R(p\wedge E)$
{\it and} $\R(p\wedge E)=\R(E,E);$\\
(III) {\it  if a Lie algebra $\g^{\h,\varphi}_3$ exists, then

 $R(\g_3^{\h,\varphi})= \R(\ker\varphi)\oplus
\R(E,\h,\varphi)\oplus\R(p\wedge E);$}\\
 (IV) {\it if a Lie algebra $\g^{\h,\psi}_4$ exists, then

 $\R(\g_4^{\h,\psi})= \R(\ker\psi)\oplus
\R(E_1,\h,\psi)\oplus\R(p\wedge E_1).$}

{\bf Remark.}
It is known (see 
\cite{B-I}) that the holonomy algebra of a weakly
irreducible, non-irreducible locally symmetric Lorentzian manifold equals  $p\wedge E=\g_2^{\{0\}}$
but this algebra can also be the
holonomy algebra of a nonlocally symmetric Lorentzian manifold (see 
\cite{Ik}).

{\bf Corollary 1.} {\it Let $\g\subset\so(V)_{\Real p}$ be a
weakly-irreducible subalgebra. Then $\g$ is a Berger algebra if and
only if $\h^\g$ is a weak-Berger algebra.}

{\bf Corollary 2.} {\it Let $\g\subset\so(V)_{\Real p}$ be a weakly-irreducible
subalgebra such that $\h^\g$ is the holonomy algebra of a Riemannian
manifold. Then $\g$ is a Berger algebra.
}

{\bf Theorem 4.}   {\it Let $\h\subset\so(E)$ be a weak-Berger
algebra. Then}\\
(I) {\it There exists an orthogonal decomposition
$E=E_0\oplus E_1\oplus\cdots\oplus E_r$ and the corresponding
decomposition into the direct sum of ideals $\h=\{0\}\oplus\h_1\oplus\h_2
\oplus\cdots\oplus\h_r$ such that $\h_i(E_j)=0$ if $i\neq j,$
$\h_i\subset\so(E_i),$ and $\h_i$ acts irreducibly on $E_i.$}\\
(II) {\it We have a decomposition
$$\P(\h)=\P(\h_1)\oplus\cdots\oplus\P(\h_r).$$
In particular,
$\h_i$ is a weak-Berger algebra for $i=1,...,k$.}

Corollary 1 and theorem 4 reduce the
classification problem for the weakly-irreducible,
non-irreducible  holonomy algebras of Lorentzian manifolds to
the classification of irreducible weak-Berger algebras.

In section 4 we will obtain the following theorem.

{\bf Theorem 5.}
{\it Let $\dim V\leq 11$, let $\g\subset\so(V)_{\Real p}$
be a weakly-irreducible subalgebra. Then
$\g$ is a Berger algebra if and only if
$\h^\g$
is the holonomy algebra of a
Riemannian manifold.
}

\section{Proof of the results}
\subsection{Proof of part (II) of theorem 3}

Let $R\in\R(\g_2^\h)$.
Above we saw that $R:\,V\wedge V\to \g_2^\h\subset V\wedge V$
is a $\eta\wedge\eta$-symmetric
linear map. It is clear that $V\wedge V=E\wedge E+p\wedge E+
q\wedge E+\Real p\wedge q$.

{\bf Lemma 1.} $R|_{(p\wedge E+\Real p\wedge q)}=0.$

Proof. By \eqref{property}, we have
$\eta(R(p\wedge v)z,w)=\eta(R(z\wedge w)p,v)$ for all $v,z,w\in V$.
Since $R(w\wedge z)\in\h$, we obtain $R(w\wedge z)p=0$ and
$\eta(R(p\wedge v)z,w)=0$. Since $\eta$ is nondegenerate,
we have $R(p\wedge v)z=0$. Thus, $R(p\wedge v)=0$. $\Box$

We see that $R$ is a linear map from $E\wedge E+q\wedge E$ to
$\h+p\wedge E$.

Let $\Gamma\subset V\wedge V$ be one of the subspaces
$E\wedge E$, $p\wedge E$, $q\wedge E$ and $\Real p\wedge q$.
Denote by
$p_\Gamma$ the projection of $V\wedge V$ onto $\Gamma$ with respect to the decomposition
$V\wedge V=E\wedge E+p\wedge E+q\wedge E+\Real p\wedge q$.

Denote by $\h^{\bot}$ the orthogonal complement of
$\h$ in $E\wedge E$,
$\h^{\bot}=\{\xi\in E\wedge
E:\eta\wedge\eta(\xi,\h)=\{0\}\}.$
Denote by
$p_\h$ the projection of $V\wedge V$ onto $\h$ with respect
to the decomposition
$V\wedge V=\h+\h^{\bot}+p\wedge E+q\wedge E+
\Real p\wedge q$.

We associate with $R$ the following linear maps:
$$R_{\h}=p_{E \wedge E} \circ R|_{ E \wedge E}:
E \wedge E \to \h,$$
$$R^{\h}_{E}=p_{p\wedge E} \circ R|_{ E \wedge E}:
E \wedge E \to p\wedge E,$$
$$R^{E}_{\h}=p_{E \wedge E} \circ R|_{q\wedge E}:
q\wedge E \to \h,$$
$$R_{E}=p_{p\wedge E} \circ R|_{q\wedge E}:
p\wedge E \to q\wedge E$$
and extend these maps to $V\wedge V$ mapping the natural complement to zero.
Then $R=R_{\h}+R^{\h}_E+
R^E_{\h}+R_E$.

{\bf Lemma 2.} $R^{\h}_{E}|_{\h^{\bot}}=0$;
$R_{\h}|_{\h^{\bot}}=0$.

Proof. Let $\theta\in \h^{\bot}$ and $\xi\in\h$.
Since the linear map $R$ is $\eta\wedge\eta$-symmetric, we have
$\eta\wedge\eta(R(\theta),\xi)=\eta\wedge\eta(R(\xi),\theta)$.
Hence,
$\eta\wedge\eta(R_{\h}(\theta)+R_{E}^{\h}(\theta),\xi)=
\eta\wedge\eta(R_{\h}(\xi)+R_{E}^{\h}(\xi),\theta)$.
Since $(p\wedge E)\,\bot\,(E\wedge E)$, $\theta\bot
R_{\h}(\xi)$, and $R_{\h}(\xi)\in\h$, we obtain
$\eta\wedge\eta(R_{\h}(\theta),\xi)=0$. Since the restriction of
the form $\eta\wedge\eta$ to $E\wedge E$ is nondegenerate, we have
$R_{\h}(\theta)=0$.

Similarly, suppose $\xi\in q\wedge E$; then
$R(\xi)=R^{E}_{\h}(\xi)+R_{E}(\xi)$.
Since $(q\wedge E)\,\bot\,(E\wedge E)$ and the restriction of the form
$\eta\wedge\eta$ to $p\wedge E+q\wedge E$ is nondegenerate, we see that
$R^{\h}_E(\theta)=0$. $\Box$

We define the linear maps $Q:\h\to E$, $T:E\to E$ and $P:E\to\h$
by conditions
\begin{equation}
R^{\h}_{E}(u\wedge v)=-1/2(p\wedge Q(u\wedge v)), \label{zam1}
\end{equation}
\begin{equation}
R_{E}(q\wedge u)=p\wedge T(u)\label{zam2}
\end{equation}
and
\begin{equation}
P(u)=R^{E}_{\h}(q\wedge u) \label{zam3}
\end{equation}
for all $u,v\in E$.

{\bf Lemma 3.} $P^*=Q$, $T^*=T$.

Proof. By \eqref{property} we have
$\eta (R(u\wedge v)q,w)=\eta(R(q\wedge w)u,v)$ for all $u,v,w\in E$.
Hence, $\eta(R^{\h}_E(u\wedge v)q,w)=\eta(R^E_{\h}(q\wedge w)u,v)$.
Using \eqref{zam1} and \eqref{zam3}, we get
$$-1/2\eta(p\wedge Q(u\wedge v))q,w)=\eta(P(w)u,v).$$
Hence, $-1/2\eta(Q(u\wedge v),w)=\eta(P(w)u,v)$.
Identity \eqref{property2} implies
$$-\eta(Q(u\wedge v),w)=\eta\wedge\eta(P(w),u\wedge v).$$

We have proved the first part of the lemma.
The second part follows from the equality
$\eta (R(q\wedge u)q,v)=\eta(R(q\wedge v)q,u)$ for all $u,v\in E$.
$\Box$

For the tensor $R$ we must check the Bianchi identity
$R(u\wedge v)w+R(v\wedge w)u+R(w\wedge u)v=0$ for all $u,v,w\in V$. It is sufficient
to check the Bianchi identity only for the basis vectors. If two of
the vectors $u$, $v$, $w$ are equal or one of the vectors $u$, $v$, $w$ equals $p$,
the identity holds trivially. Thus it is sufficient to check
the Bianchi identity in the two cases: $u,v,w\in E$; $u,v\in E$, $w=q$.
We do this in the following lemma.

{\bf Lemma 4.} $R_{\h}\in\R(\h),$ $P\in\P(\h)$.

Proof. Let us write the Bianchi identity for $u,v,w\in E$;
$R(u\wedge v)w+R(v\wedge w)u+R(w\wedge u)v=0$.
From the equalities like
$R(u\wedge v)=R_{\h}(u,v)+R^{\h}_E(u,v)$ and
\eqref{zam1} it follows that
$R_{\h}(u\wedge v)w+R_{\h}(v\wedge w)u+
R_{\h}(w\wedge u)v-
1/2(p\wedge Q(u\wedge v))w-1/2(p\wedge Q(v\wedge w))u-1/2(p\wedge Q(w\wedge u))v=0$.
Since $(p\wedge Q(u\wedge v))w=-\eta(Q(u\wedge v),w)p\in\Real p$ and
$R_{\h}(u\wedge v)w\in E$, we obtain
\begin{equation}
R_{\h}(u\wedge v)w+R_{\h}(v\wedge w)u+R_{\h}(w\wedge u)v=0
\label{B1}\end{equation} and
\begin{equation}
\eta(Q(u\wedge v),w)+\eta(Q(v\wedge w),u)+\eta(Q(w\wedge u),v)=0.\label{B2}\end{equation}
Identity \eqref{B1} shows that
$R_{\h}\in\R(\h)$.

Now write the Bianchi identity for $u,v\in E$ and $q$;
$R(u\wedge v)q+R(v\wedge q)u+R(q\wedge u)v=0$.
Hence, $R^{\h}_E(u\wedge v)q+R^E_{\h}(v\wedge q)u+
R^E_{\h}(q\wedge u)v=0$.
Combining this with \eqref{zam1} and \eqref{zam3}, we obtain
$-1/2(p\wedge Q(u\wedge v))q-P(v)u+P(u)v=0$. Hence, $1/2Q(u\wedge v)q+P(v)u-P(u)v=0$.
Combining this with \eqref{B2} and using the fact that for $z\in E$
the endomorphism $P(z)$ is $\eta$-skew symmetric, we obtain
\begin{equation}
\eta(P(u)v,w)+\eta(P(v)w,u)+\eta(P(w)u,v)=0 \text{ }\text{\it for all } u,v,w\in E.
\label{Pprop}\end{equation}
Identity \eqref{Pprop} implies $P\in\P(\h)$.
 $\Box$

We put $R^P=R_E^\h+R_\h^E$, $R^T=R_E$.
Then $R^P\in\R(E,\h),$ $R^T\in\R(p\wedge E)$, and
$R=R_\h+R^P+R^T$.
We have proved that
$\R(\g_2^\h)\subset\R(\h)\oplus
\R(E,\h)\oplus\R(p\wedge E).$
Now we check the inverse inclusion.

Suppose $T\in S^2(E)$, i.e. $T:E\to
E$ is a linear map such that $T^*=T$. Put $R^T(q\wedge u)=p\wedge T(u)$,
$R^T(p\wedge q)=R^T(p\wedge u)=R^T(u\wedge v)=0$ for all $u,v\in E$.
Then $R^T\in\R(p\wedge E)\subset\R(\g_2^\h)$.

Let $P\in\P(\h)$ and $R^P\in\R(E,\h)$.
We must check the Bianchi identity for the
tensor $R^P$.
For $u,v,w\in E$ the identity follows from \eqref{Pprop}.

Suppose $u,v\in E$; then $R^P(u\wedge v)q+R^P(v\wedge q)u+R^P(q\wedge u)v=$
$-1/2(p\wedge P^*(u\wedge v))q-P(v)u+P(u)v=-1/2P^*(u\wedge v)-P(v)u+P(u)v.$

Suppose $w\in E$; then $\eta(-1/2P^*(u\wedge v)-P(v)u+P(u)v,w)=
\eta(-1/2P^*(u\wedge v),w)-\eta(P(v)u,w)+\eta(P(u)v,w)=$
$-1/2\eta(P^*(u\wedge v),w)+\eta(P(v)w,u)+\eta(P(u)v,w)=$
$-1/2\eta(P^*(u\wedge v),w)-\eta(P(w)u,v)=$ $-1/2(\eta(P^*(u\wedge
v),w))+\eta\wedge\eta(P(w),u\wedge v)=0.$ Above we used the fact
that $P(v)\in\h\subset\so(E)$, \eqref{property2} and \eqref{property}. Since
the restriction of $\eta$ to $E$ is nondegenerate, we obtain
$R^P(u\wedge v)q+R^P(v\wedge q)u+R^P(q\wedge u)v=0$.

Now we have
$\R(\h)\oplus\R(E,\h)\oplus\R(p\wedge E)
\subset\R(\g_2^\h).$
Thus,
$\R(\g_2^\h)=
\R(\h)\oplus\R(E,\h)\oplus\R(p\wedge E).$

\subsection{Proof of part (I) of theorem 3}

Let $R\in\R(\g_1^\h)$.
Similarly to lemma 1, we can prove that $R|_{p\wedge
E}=0.$ Hence $R$ is a linear map from
$E\wedge E+q\wedge E+\Real p\wedge q$ to
$\h+p\wedge E+\Real p\wedge q.$

We can define the linear maps
$R_{\h},\,R^{\h}_E,\,R^E_{\h}$ and $R_E$
as in section 3.1. It is easily shown that the map
$R_2=R_{\h}+R^{\h}_E+R^E_{\h}+R_E$
is a curvature tensor of type the Lie algebra $\g_2^\h$.

We define the following linear maps:

$$ R^{\h}_\Real=p_{\Real p \wedge q} \circ R|_{ E
\wedge E}: E \wedge E \to \Real p \wedge q,$$
$$ R^\Real_{\h}=p_{E \wedge E} \circ
R|_{ \Real p \wedge q}: \Real p \wedge q \to E \wedge E,$$
$$R^\Real_{E}=p_{ p \wedge E } \circ R|_{ \Real p \wedge q }: \Real p \wedge
q \to p \wedge E , $$
$$ R^E_\Real=p_{ \Real p \wedge q} \circ R|_{
q \wedge E}: q \wedge E \to \Real p \wedge q, $$
$$R_\Real=p_{\Real
p \wedge q} \circ R|_{\Real p \wedge q}: \Real p \wedge q\to
\Real p \wedge q$$

and extend these maps to $V\wedge V$ sending the natural complement to zero.

We have $R=R_2+R^{\h}_\Real+
R^\Real_{\h}+R^\Real_{E}+R^E_\Real+R_\Real$

{\bf Lemma 5.} $R^{\h}_\Real=0,$ $R^\Real_{\h}=0.$

Proof. Let as write the Bianchi identity for vectors $u,v\in E$
and $p$; $R(u\wedge v)p+R(v,p)u+R(p\wedge u)v=0.$ Since $R(v\wedge p)=R(p\wedge u)
=0$ and
$R(u\wedge v)p=R^{\h}_\Real(u\wedge v)p,$ we see that
$R^{\h}_\Real(u\wedge v)=0$.

Using \eqref{property}, we get $\eta(R(p\wedge q)u,v)=\eta(R(u\wedge v)p,q)=0.$
Since $u$ and $v$ are arbitrary and the restriction of $\eta$ to $E$
is not degenerate, we obtain $R^\Real_{\h}=R(p\wedge q)|_E=0.$
$\Box$

Define the linear maps $L:E\to\Real$ and $K:\Real\to E$ by
conditions
\begin{equation}
R^E_\Real(q\wedge u)=L(u)p\wedge q \text{ for all } u\in E, \label{zam4}\end{equation}
\begin{equation}
R^\Real_E(a\,p\wedge q)=p\wedge K(a) \text{ for all } a\in\Real. \label{zam5}\end{equation}

{\bf Lemma 6.} $K=L^*.$

The proof is similar to the proof of lemma 3. $\Box$

Let $\lambda$ be a real number such that
$R_\Real(p\wedge q)=\lambda p\wedge q$. Put
$R^L=R^E_\Real+R_E^\Real$ and $R^\lambda=R_\Real$. We see that
$R^L\in\R(E,\Real)$, $R^\lambda\in\R(\Real,\Real)$, and
$R=R_2+R^L+
R^\lambda\in\R(\g_2^\h)\oplus \R(E,\Real)\oplus\R(\Real,\Real)$.
Thus we have $\R(\g_1^\h)\subset\R(\g_2^\h)\oplus \R(E,\Real)\oplus\R(\Real,\Real)$.
The inverse inclusion is obvious.

\subsection{Proof of part (III) of theorem 3}

 We have $\g_2^\h\subset\g_3^{\h,\varphi}\subset\g_1^h.$
Suppose $R\in\R(\g_3^{\h,\varphi}),$ then we have
$R\in\R(\g_1^\h).$ From section 3.2 it follows that
$R=R_2+R^E_\Real+R^\Real_E+R_\Real,$ where $R_2=R_\h+
R^\h_E+R^E_\h+R_E\in\R(\g_2^\h).$

We claim that $R_\Real=0.$ Indeed, there exists a $\lambda\in\Real$ such
that $R_\Real(p\wedge q)=\lambda\,p\wedge q;$ we have
$R(p\wedge q)=R_\Real(p\wedge q)+R^\Real_E(p\wedge q)=\lambda\,p\wedge
q+R^\Real_E(p\wedge q).$ Since $R^\Real_E(p\wedge q)\in p\wedge E\subset\g_3^{\h,\varphi}$
and
$\Real p\wedge q\cap \g_3^{\h,\varphi}=\{0\}$, we see that $\lambda=0.$

For any $u\in E$ we have
$R(q\wedge u)=R^E_\h(q\wedge u)+R^E_\Real(q\wedge u)\in \h+\Real p\wedge q.$
Since $R(q\wedge u)\in \g_3^{\h,\varphi},$ we see that
$R^E_\Real(q\wedge u)=\varphi(R^E_\h(q\wedge u))p\wedge q.$
Hence, $R^E_\Real=p\wedge q\varphi\circ R^E_\h.$
Using \eqref{zam3} and \eqref{zam4}, we get $L(u)=\varphi(P(u)).$
Hence for any $a\in\Real$ we have
$L^*(a)=P^*(\varphi^*(a)).$
Using this, \eqref{zam5} and lemma 6, we obtain
$R^\Real_E(a\,p\wedge q)=p\wedge L^*(a)
=p\wedge P^*(\varphi^*(a)).$
Using \eqref{zam1}, we obtain
$R^\Real_E(a\,p\wedge q)=R^\h_E\circ\varphi^*(a)$ for all
$a\in\Real$.

Suppose $u,v\in E;$ then we have
$R(u\wedge v)=R_\h(u\wedge v)+R^\h_E(u\wedge v)\in\g_3^{\h,\varphi}.$
Hence, $R_\h(u\wedge v)\in\ker\varphi.$ We see that
$R_\h\in\R(\ker\varphi).$ Put
$R_{\ker\varphi}=R_\h.$

Put $R^P=R^E_{\h}+R^{\h}_E+R^E_{\Real}+R^{\Real}_E$.
We see that $R^P\in\R(E,\h,\varphi)$.
Thus, $R=R_{\ker\varphi}+R^P+R_E\in\R(\ker\varphi)\oplus
\R(E,\h,\varphi)\oplus\R(p\wedge E).$

Conversely, let $R=R_{\ker\varphi}+R^P+R_E\in\R(\ker\varphi)\oplus
\R(E,\h,\varphi)\oplus\R(p\wedge E).$
From section 3.2 it follows that
$R\in\R(\g_1^\h).$
Since for any $u,v\in V$ we have $R(u\wedge v)\in\g_3^{\h,\varphi},$ we see
that $R\in\R(\g_3^{\h,\varphi}).$

\subsection{Proof of part (IV) of theorem 3}

By definition, $\g_2=\h+p\wedge E_1$
is a Lie algebra of type 2.
We have
$\g_2\subset\g_4^{\h,\psi}\subset\g_2^\h.$
Suppose that $R\in\R(\g_4^{\h,\psi})$, then we have
$R\in\R(\g_2^\h).$ From section 3.1 it
follows that $R=R_{\h}+R^{E}_{\h}+R^{\h}_{E}+R_{E}$.
There exists a $P\in\P(\h)$ such that
$R^{E}_{\h}(q\wedge u)=P(u)$ for all $u\in E$.
Let $u_1,v_1\in E_1$, $u_0\in E_0$.
We have
$\eta(P(u_1)v_1,u_0)+\eta(P(v_1)u_0,u_1)+\eta(P(u_0)u_1,v_1)=0.$
Since $\h(E_0)=\{0\}$ and $\h(E_1)\subset E_1$, we see that
$\eta(P(u_0)u_1,v_1)=0$ for all $u_1,v_1\in E_1$, $u_0\in E_0$.
Hence, $P(E_0)=\{0\}$, $P^*(\h)\subset E_1$.
We can write
$R^{E}_{\h}=R^{E_1}_{\h}$ and
$R^{\h}_{E}=R^{\h}_{E_1}$.

We consider the following linear maps:
$$R_{E_1}=p_{p\wedge E_1}\circ R_E|_{q\wedge E_1}:q\wedge E_1\to p\wedge E_1 ,$$
$$R^{E_0}_{E_1}=p_{p\wedge E_1}\circ R_E|_{q\wedge E_0}:q\wedge E_0\to p\wedge E_1,$$
$$R^{E_1}_{E_0}=p_{p\wedge E_0}\circ R_E|_{q\wedge E_1}:q\wedge E_1\to p\wedge E_0,$$
$$R_{E_0}=p_{p\wedge E_0}\circ R_E|_{q\wedge E_0}:q\wedge E_0\to p\wedge E_0$$

and extend these maps to $V\wedge V$ mapping the complementary
subspace to zero.
Obviously, $R_{E}=R_{E_1}+R^{E_0}_{E_1}+R^{E_1}_{E_0}+R_{E_0}$.

We claim that $R_{E_0}=0.$ Indeed, for $u_0\in E_0$ we have
$R(q\wedge u_0)=R_{E}(q,u_0)=R^{E_0}_{E_1}(q\wedge
u_0)+R_{E_0}(q\wedge u_0)\in\g_4^{\h,\psi};$ since
$R^{E_0}_{E_1}(q\wedge u_0)\in p\wedge E_1\subset\g_4^{\h,\psi}$,
$R_{E_0}(q\wedge u_0)\in p\wedge E_0$, and
$p\wedge E_0\cap\g_4^{\h,\psi}=\{0\}$,
we obtain $R_{E_0}(q\wedge u_0)=0,$ hence, $R_{E_0}=0.$

As in section 3.3, we can prove that
$R_{\h}\in\R(\ker\psi),$
$R^{E_1}_{E_0}=p\wedge\psi\circ R^{E_1}_{\h}$, and
$R^{E_0}_{E_1}(q\wedge u_0)=R^{\h}_{E_1}\circ\psi^*(u_0)$ for all
$u_0\in E_0$.

Put $R^P= R^{E_1}_{\h}+R^{\h}_{E_1}+R^{E_1}_{E_0}+R^{E_0}_{E_1}$.
Thus we have
$R=R_{\ker\psi}+R^P+R_{E_1}
\in\R(\ker\psi)\oplus \R(E_1,\h,\psi)\oplus\R(p\wedge E_1).$

Conversely, let
$R=R_{\ker\psi}+R^P+R_{E_1}
\in\R(\ker\psi)\oplus \R(E_1,\h,\psi)\oplus\R(p\wedge E_1).$
From section 3.1 it follows that
$R\in\R(\g_2^\h).$
Since for any $u,v\in V$ we have $R(u\wedge v)\in\g_4^\h,$ we see
that $R\in\R(\g_4^\h).$

This concludes the proof of theorem 3. $\Box$

\subsection{Proof of theorem 4}

(I) Suppose $\h$ preserves a proper subspace $E_1\subset E$, then
$\h$ preserves the orthogonal complement $E_1^\bot$ to $E_1$ in
$E$. Put $E_2=E_1^\bot$. We have $E=E_1\oplus E_2$,
$\h(E_1)\subset E_1$, and $\h(E_2)\subset E_2$. Put
$\h_1=\{\xi\in\h: \xi(E_2)=\{0\}\}$ and
$\h_2=\{\xi\in\h:\xi(E_1)=\{0\}\}$. Obviously, $\h_1$ and $\h_2$
are ideals in $\h$ and $\h_1\cap\h_2=\{0\}$.

Let $P\in\P(\h)$.
Let $u_1,v_1\in E_1$, $u_2\in E_2$.
We have $\eta(P(u_1)v_1,u_2)+\eta(P(v_1)u_2,u_1)+\eta(P(u_2)u_1,v_1)=0.$
Since $\h(E_1)\subset E_1$ and $\h(E_2)\subset E_2$, we have
$\eta(P(u_2)u_1,v_1)=0$ for all $u_1,v_1\in E_1$, $u_2\in E_2$.
We see that $P(E_1)\subset \h_1$ and $P(E_2)\subset \h_2$.
Hence, $L(\P(\h))\subset\h_1\oplus\h_2.$
Combining this with the equality $L(\P(\h))=\h$,
we obtain $\h=\h_1\oplus\h_2$.

(II) Let $P\in\P(\h)$. As above, we can prove that
$P(E_1)\subset \h_1$ and $P(E_2)\subset \h_2$.
By definition, put $P_1=P|_{E_1},$ $P_2=P|_{E_2}.$
It is clearly that $P_1\in\P(\h_1)$, $P_2\in\P(\h_2)$, and
$P=P_1+P_2$.

Conversely, for any $P_1\in\P(\h_1)$ and $P_2\in\P(\h_2)$ we have
$P=P_1+P_2\in\P(\h).$

The proof of the theorem follows easily by induction over the
number of the summands. $\Box$

\subsection{Proof of corollaries}

Let $\h\subset\so(E)$ be a subalgebra and $R\in\R(\h)$.
We claim that for any $z\in E$ the tensor $P$ defined by
$P(\cdot)=R(\cdot\wedge z)$ belongs to $\P(\h)$.
Indeed, we have $R(u\wedge v)w+R(v\wedge w)u+R(w\wedge u)v=0$
for all $u,v,w\in E.$ Multiplying both sides innerly by $z\in E,$
we obtain
$\eta(R(u\wedge v)w,z)+\eta(
R(v\wedge w)u,z)+\eta(R(w\wedge u)v,z)=0$ for
all $u,v,w,z\in E.$
Using \eqref{property}, we get
$\eta(R(w\wedge z)u,v)+\eta(
R(u\wedge z)v,w)+\eta(R(v\wedge z)w,u)=0.$

From theorem 3 and the claim it
follows that for any weakly-irreducible algebra
$\g$ we have $p_{\so(E)}(L(\R(\g)))=L(\P(\h^\g))$.

Corollary 1 follows from  the following
obvious facts: $L(\R(p\wedge E))=p\wedge E$ and
$L(\R(\Real,\Real))=\Real p\wedge q$.

Corollary 2 follows from corollary 1 and the above claim. $\Box$

\section{Examples}

Above we have reduced the classification problem for
Berger algebras of Lorentzian manifolds to the classification of
irreducible weak-Berger algebras.

Suppose we have the full list of  irreducible weak-Berger algebras.
Corollary 1 and theorem 4 imply that the full list of
Berger algebras of Lorentzian manifolds
 can be obtained in the following way.

For each Euclidean space $E$ we must consider all orthogonal
decompositions $E=E_0\oplus E_1\oplus\cdots\oplus E_r$ such that
$2\leq\dim E_1\leq\cdots\leq\dim E_r$, and for each
Euclidean space $E_i$ all the irreducible weak-Berger algebras
$\h_i\subset\so(E_i)$.
From theorem 4 it follows that the algebras
$\h=\h_1\oplus\cdots\oplus\h_r$
exhaust all weak-Berger algebras. Corollary 1 implies that
the Lie algebras $\g_1^\h$, $\g_2^\h$, $\g_3^{\h,\varphi}$ and
$\g_4^{\h,\psi}$ (if $\g_3^{\h,\varphi}$ and
$\g_4^{\h,\psi}$ exist) exhaust all Berger algebras.

Below we list all the irreducible subalgebras
$\h\subset\so_n$ for $n\leq 9$
and state the result of computing of the spaces $\P(\h)$ for
algebras that are not the holonomy algebras of Riemannian
manifolds.

Since $\h\subset\so(E)$, the Lie
algebra $\h$ is compact. Hence,
$\h=\h'\oplus\mathfrak{z}(\h),$
where $\h'$ is a compact semisimple ideal,
$\mathfrak{z}(\h)$ is an Abelian ideal.
Since the subalgebra $\h\subset\so(E)$ is
irreducible, by Schur lemma the center $\z(\h)$ has dimension $0$
or $1$.

It is known that if a subalgebra
$\h\subset\so(E)$ is irreducible, the
subalgebra $\h'\subset\so(E)$ is irreducible too (see 
\cite{Z-S}).
Let $\h\subset\so(E)$
be a  semisimple irreducible subalgebra. Denote by
$\mathfrak{z}_{\so(E)}(\h)$ the centralizer of
$\h$ in $\so(E).$
If $\mathfrak{z}_{\so(E)}(\h)\neq\{0\},$ then
for each one-dimensional subspace
$\t\subset\mathfrak{z}_{\so(E)}(\h)$
the Lie algebra $\h\oplus\t$ is a compact non-semisimple
irreducible subalgebra of $\so(E).$ Hence it is
sufficient to get the list of all semisimple irreducible
subalgebras of $\so(E).$

The classification of irreducible representations of compact
semisimple Lie algebras is well known, see for example 
\cite{V-O}.
Any irreducible real representation  $\pi:\h\to\mathfrak{gl}(E)$ of a
real semisimple Lie algebra $\h$ can be obtained in the
following way.

Suppose we have a complex irreducible representation
$\rho:\mathfrak{f}\to\mathfrak{gl}(U)$
of a complex semisimple Lie algebra $\mathfrak{f}$.
 Let $\h\subset\mathfrak{f}$ be a compact real form of $\mathfrak{f}$
($\h$ is unique up to conjugation).

There are the following three cases:

1) The representation $\rho$ is self-dual and orthogonal
(i.e. $\rho\thicksim\rho^*,$ where
$\rho^*:\mathfrak{f}\to\mathfrak{gl}(U^*)$ is the dual
representation for $\rho$ and $\rho$ admits an invariant not degenerate
symmetric bilinear form)

2) The representation $\rho$ is self-dual and symplectic
(i.e. $\rho$ admits an invariant not degenerate
skew symmetric bilinear form)

3) The representation $\rho$ is not self-dual.

The first condition holds if and only if the representation $\rho$
admits a real form $J,$ i.e. $J:U\to U$ is a $\Real$-linear map
such that $J(iu)=-iJ(u)$ for all $u\in U,$
$J^2=\text{id},$ and $J\rho(\xi)=\rho(\xi)J$ for all $\xi\in\mathfrak{f}.$
In this case the real representation $\pi=\rho|_{\h}$ in the
realificated vector space $U^{\Real}$ preserves the space $U^J=\{u\in
U:\,J(u)=u\}$ and acts irreducibly on $U^J.$ We get an irreducible
real representation $\pi:\h\to\mathfrak{gl}(U^J).$

In the cases 2) and 3) the real representation
$\pi=\rho|_{\h}:\h\to\mathfrak{gl}(U^\Real)$
is irreducible.
In the case 2) we have $\pi(\h)\subset\sp_m,$
in the case 3) we have
$\pi(\h)\subset\su_n,$ where $2n=4m=\dim_{C} U.$

Let $\pi:\h\to\mathfrak{gl}(E)$ be an irreducible real
representation of a compact Lie algebra $\h$ in a real
vector space $E.$ Since $\h$ is compact, we see that $\pi$
admits an invariant symmetric positively definite bilinear form.
This form is unique up to non-zero real factor. The linear space
$E$ is an Euclidean space with respect to this form and we can
write $\pi(\h)\subset\so(E).$

Any irreducible complex representation of a complex semisimple Lie algebra
$\mathfrak{f}$ is uniquely defined (up to equivalence of representations) by
its highest weight $\Lambda.$ The highest weight $\Lambda$ can be
given by the labels $\Lambda_1,...,\Lambda_l$ on the Dynkin
diagram of the Lie algebra $\mathfrak{f}$ ($l=\rk(\mathfrak{f}$)).
There exists a criteria for a complex representation
to be orthogonal, symplectic or self-dual in terms of the highest
weight.

In the table 1 we list all irreducible  subalgebras
$\pi(\h)\subset\so_n$ for $n\leq 9$.

The second column of the table contains the irreducible holonomy algebras of
Riemannian manifolds.  The third column of the table
contains algebras that are not the holonomy algebras of
Riemannian manifolds.

Let $\h$ be a compact semisimple Lie algebra. We denote
by $\pi^\mathbf{K}_{\Lambda_1,...,\Lambda_l}:\h\to\Real^n$
the real irreducible representation  that
is obtained as above from the complex representation
$\rho_{\Lambda_1,...,\Lambda_l}:\mathfrak{h}(\mathbf{C})
\to \mathfrak{gl}(U)$ with the highest
weight $\Lambda,$ here $\Lambda_1,...,\Lambda_l$ are the labels
of $\Lambda$ and $\mathbf{K}=\Real,\,\mathbf{C}$ or $\mathbf{H}$ if
the representation $\rho_{\Lambda_1,...,\Lambda_l}$ is orthogonal,
not self-dual or symplectic respectively.
If the representation $\rho_{\Lambda_1,...,\Lambda_l}$
is orthogonal, then $n=\dim_{\mathbf{C}}U,$ otherwise we have
$n=2\dim_{\mathbf{C}}U.$
We denote by $\t$ the 1-dimensional Lie algebra $\Real$.

\newpage
{\bf Table 1.} Irreducible subalgebras of $\so_n$

\vskip 0.3cm
\begin{tabular}{|c|c|c|}
\hline

  $n$ & irreducible holonomy algebras of $n$-dimensional &
 other irreducible \\

 & Riemannian manifolds & subalgebras of $\so_n$ \\ \hline

  $n=1$ &   &  \\ \hline
  $n=2$ & $\so_2$    &  \\ \hline
  $n=3$ & $\pi^\Real_2(\so_3)$   &  \\ \hline
  $n=4$ &$\pi^\Real_{1,1}(\so_3+\so_3)$,
  $\pi^\mathbf{C}_1(\su_2)$,
  $\pi^\mathbf{C}_1(\su_2)\oplus\t$

    &  \\ \hline
  $n=5$ & $\pi^\Real_{1,0}(\so_5)$,
  $\pi^\Real_4(\so_3)$ & \\ \hline
  $n=6$ & $\pi^\Real_{1,0,0}(\so_6)$,
  $\pi^\mathbf{C}_{1,0}(\su_3)$,
   $\pi^\mathbf{C}_{1,0}(\su_3)\oplus\t$
    &  \\ \hline
  $n=7$ & $\pi^\Real_{1,0,0}(\so_7)$, $\pi^\Real_{1,0}(\g_2)$
    & $\pi^\Real_6(\so_3)$   \\ \hline

   $n=8$ & $\pi^\Real_{1,0,0,0}(\so_8)$,
    $\pi^\mathbf{C}_{1,0}(\su_4)$,
    $\pi^\mathbf{C}_{1,0}(\su_4)\oplus\t$,
  $\pi^\mathbf{H}_{1,0}(\sp_2)$,
   
&
   $\pi^\mathbf{C}_3(\so_3)$,
    $\pi^\mathbf{C}_3(\so_3)\oplus\t$

     \\
 & $\pi^\Real_{1,0,1}(\sp_2\oplus\sp_1),$
    $\pi^\Real_{1,3}(\so_3\oplus\so_3)$,
  $\pi^\Real_{1,1}(\su_3)$

    & $\pi^\mathbf{H}_{1,0}(\sp_2)\oplus\t$ \\      \hline

$n=9$ & $\pi^\Real_{1,0,0,0}(\so_9)$,
$\pi^\Real_{2,2}(\so_3\oplus\so_3)$&
$\pi^\Real_8(\so_3)$   \\ \hline
  \end{tabular}

\vskip 0.3cm

Now we must verify the equality $L(\P(\h))=\h$
for algebras from the third column of the table.

Let $\h\subset\so(E)$ be a subalgebra.
We claim that $L(\P(\h))$ is an ideal in $\h$.
Indeed, let $P\in\P(\h)$ and $\xi\in\h$;
put $P_\xi(u)=-\xi\circ P(u)+P(u)\circ\xi+P(\xi u)$ for all $u\in E$.
It can be easily checked that $P_\xi\in\P(\h)$.
We see that $[P(u),\xi]=P_\xi(u)-P(\xi u)$ for all $u\in E$, $\xi\in\h$.

Suppose $\mathfrak{h}\subset\mathfrak{so}(E)$ is an irreducible
subalgebra. We have
$\mathfrak{h}=\mathfrak{h}'\oplus\mathfrak{z}(\mathfrak{h})$ and
$\dim\mathfrak{z}(\mathfrak{h})=0$ or $1$. Since $\mathfrak{h}'$ is
semisimple, we have
$\mathfrak{h}'=\mathfrak{h}_1\oplus\ldots\oplus\mathfrak{h}_r$,
where $\mathfrak{h}_i\subset \mathfrak{h}'$ are simple ideals. Any
ideal $\tilde{\mathfrak{h}}\subset\mathfrak{h}$ is a sum of some
of the ideals $\mathfrak{h}_1, \ldots, \mathfrak{h}_r$ and, may
be,
$\mathfrak{z}(\mathfrak{h})$. From the above claim and the obvious
equality
$\mathcal{P}(\mathfrak{h})=\mathcal{P}(L(\mathcal{P}(\mathfrak{h})))$
it follows that $L(\P(\h))=\h$ if and only if
$\mathcal{P}(\mathfrak{h})\neq\{0\}$ and there is no proper ideal
$\tilde{\mathfrak{h}}\subset\mathfrak{h}$ such that
$\mathcal{P}(\mathfrak{h})=\mathcal{P}(\tilde{\mathfrak{h}})$.} In
particular, if $\h$ is simple, then $L(\P(\h))=\h$ if and only if
$\mathcal{P}(\mathfrak{h})\neq\{0\}$.

Let $\h\subset\so(E)$ be a subalgebra and $N=\dim\h$. Denote by
$C^1,...,C^N$ a basis of the vector space $\h$.
Let $P\in\text{\rm Hom}(E,\h)$ and $P_{\alpha\,i}$
($\alpha=1,...,N;
i=1,...,n$) be real numbers such that $P(e_i)=\sum_{\alpha=1}^N
P_{\alpha\,i}C^\alpha$ for $i=1,...,n$.
We have $P\in\P(\h)$ if and only if
$$\sum_{\alpha=1}^N(
P_{\alpha\,i}C^\alpha_{j\,k}+P_{\alpha\,j}C^\alpha_{k\,i}+P_{\alpha\,k}C^\alpha_{i\,j})=0$$
for all $i,j,k$ such that $1\leq i<j<k\leq n$
($\{C^\alpha_{i\,j}\}_{i,j=1}^n$ is the matrix of the endomorphism $C^\alpha$).

Thus the space $\P(\h)$ can be found as the solution of the system
of $n(n-1)(n-2)/6$ equations in $nN$ unknowns.
We used computer program Mathematica 4.1 to solve such systems for
algebras from the third column of  table 1.

The result is\\ $\mathcal{P}(\pi^\Real_6(\so_3))=\{0\},$
$\mathcal{P}(\pi^\mathbf{C}_3(\so_3))=\{0\},$
$\mathcal{P}(\pi^\mathbf{C}_3(\so_3)\oplus\t)=\{0\},$
$\dim(\mathcal{P}(\pi^\mathbf{H}_{1,0}(\sp_2)\oplus\t))=40$,
$\mathcal{P}(\pi^\Real_8(\so_3))=\{0\}$. We also have
$\dim(\mathcal{P}(\pi^\mathbf{H}_{1,0}(\sp_2)))=40$.

Since $\mathcal{P}(\pi^\mathbf{H}_{1,0}(\sp_2))=
\mathcal{P}(\pi^\mathbf{H}_{1,0}(\sp_2)\oplus\t),$ we see
that
$L(\mathcal{P}(\pi^\mathbf{H}_{1,0}(\sp_2)\oplus\t))=
\pi^\mathbf{H}_{1,0}(\sp_2)$.
Hence the algebra
$\pi^\mathbf{H}_{1,0}(\sp_2)\oplus\t$ is
not a weak-Berger algebra.

In table 2  we list all irreducible
weak-Berger algebras $\h\subset\so_n$ $(n\leq 9)$.
This list coincides with the list of the irreducible holonomy
algebras of Riemannian manifolds. We use the standard notation.
In the table $\otimes$
stands for the tensor product of representations;
$\underline{\otimes}$
stands for the highest irreducible component of the corresponding
product.

\vskip 0.3cm

{\bf Table 2.} Irreducible weak-Berger algebras

\vskip 0.3cm

\begin{tabular}{|c|c|}
\hline
 n & Irreducible weak-Berger subalgebras of $\so_n$\\ \hline
  n=1 &    \\ \hline
  n=2 & $\so_2$   \\ \hline
  n=3 & $\so_3$  \\ \hline
  n=4 & $\so_4$,
  $\su_2$,
  $\u_2$

    \\ \hline
  n=5 & $\so_5$,
  $\underline{\otimes}^2\so_3$  \\ \hline
  n=6 & $\so_6$,
  $\su_3,$
   $\u_3$
  \\ \hline
  n=7 & $\so_7$, $\g_2$
     \\ \hline

   n=8 & $\so_8$,
    $\su_4$,
    $\u_4$,

$\sp_2$,
    $\sp_2\otimes\sp_1$,
    $\underline{\otimes}^3\so_3\otimes\so_3$,

    $\underline{\otimes}^2\su_3$

    \\      \hline

n=9 & $\so_9$,
$\so_3\otimes\so_3$
 \\      \hline
  \end{tabular}

\vskip 0.4cm

Recall that the holonomy group of an indecomposable Lorentzian
manifold can be not closed. In 
\cite{B-I} it was shown that
the connected Lie subgroups of $SO_{1,n+1}$ corresponding to Lie
algebras of type 1 and 2 are closed; the connected Lie subgroup
of $SO_{1,n+1}$ corresponding to a Lie
algebra of type 3 (resp. 4) is closed if and only if the connected
Lie subgroup of $SO_n$ corresponding to the subalgebra
$\ker\varphi\subset\z(\h)$ (resp. $\ker\psi\subset\z(\h)$) is
closed. We give a criteria for Lie groups corresponding to Lie
algebras of type 3 and 4 to be closed in terms of the Lie algebras
$\ker\varphi$ and $\ker\psi$.

Let $\h\subset\so_n$ be a weak-Berger algebra such that $\z(\h)\neq\{0\}$.
Denote by $T$ the connected Lie subgroup of $SO_n$ corresponding to the Lie
subalgebra $\z(\h)\subset\so_n$. Since $\h$ is a sum of
irreducible weak-Berger algebras, we see that the Lie subgroup $T$
is closed. Hence $T$ is compact and $T$ is isomorphic to the
torus of dimension $k=\dim\z(\h)$. Thus, $T=S_1\times\cdots\times S_k$,
where $S_i$ is
a closed Lie subgroup of $T$ isomorphic to the unit circle $S^1$ ($i=1,...,k$).
Denote by $u_1,...,u_k$ the tangent vectors to $S_1,...,S_k$
respectively corresponding to a unit
tangent vector at the unity element of the circle. The vectors
$u_1,...,u_k$ form a basis of $\z(\h)$.

Let $\tilde{\t}\subset\z(\h)$ be a subalgebra and
$\tilde{T}\subset T$ the corresponding connected Lie subgroup.
We claim that {\it the subgroup $\tilde{T}$ is closed if and only if
there exists a basis $v_1,...,v_l$ ($l=\dim\tilde{\t}$) of the
vector space $\tilde{\t}$ such that the coordinates of the vector $v_i$
with respect to the basis $u_1,...,u_k$ are integer for all
$i=1,...,l$.}

For $l=1$ this statement was proved in 
\cite{V-O}.

Let $l>1$. Suppose that there exists a basis of $\tilde{\t}$
as above. Denote by $\tilde{S_1}$,...,$\tilde{S_l}$ the connected
Lie subgroups of $T$ corresponding  to the  subalgebras
$\Real v_1,...,\Real v_l\subset\z(\h)$. The Lie
subgroups  $\tilde{S_1}$,...,$\tilde{S_l}$ are closed.
Hence the Lie groups $\tilde{S_1}$,...,$\tilde{S_l}$ are compact
and isomorphic to the unit circle. Denote these isomorphisms by
$f_1$,...,$f_l$. Put $T^l=S^1\times\cdots\times S^1$.
Define a map $f:T^l\to T$ by putting
$f(x_1,...,x_l)=f_1(x_1)\cdots f_l(x_l)$, where $x_i\in S^1$.
We have $f(T^l)=\tilde{T}$, hence $\tilde{T}$ is closed in $T$.
The inverse statement is obvious.


\bibliographystyle{unsrt}

\end{document}